\documentclass[11pt]{article}

\usepackage{amsfonts,amssymb,amstext,amsmath,amsthm,latexsym,color,epsfig, tikz}
\newtheorem{theorem}{Theorem}

\newtheorem*{HJ}{Hales--Jewett theorem}

\makeatletter
\def\Ddots{\mathinner{\mkern1mu\raise\p@
\vbox{\kern7\p@\hbox{.}}\mkern2mu
\raise4\p@\hbox{.}\mkern2mu\raise7\p@\hbox{.}\mkern1mu}}
\makeatother

\title{\vspace{-0.7cm}Monochromatic combinatorial lines of length three}
\author{David Conlon\thanks{Department of Mathematics, California Institute of Technology, Pasadena, CA 91125, USA. E-mail: {\tt dconlon@caltech.edu}. Research supported by ERC Starting Grant 676632 and by NSF Award DMS-2054452.}}
\date{}
\begin{document}
\maketitle

\begin{abstract}
We show that there is a positive constant $c$ such that any colouring of the cube $[3]^n$ in $c \log \log n$ colours contains a monochromatic combinatorial line. 
\end{abstract}

\section{Introduction}

The Hales--Jewett theorem~\cite{HJ63} is a central result in Ramsey theory, an abstract version of van der Waerden's theorem~\cite{vdW27} saying that finite colourings of high-dimensional cubes contain monochromatic lines.
To state the Hales--Jewett theorem formally, we consider the cube $[m]^n$ with $[m] = \{1,2, \dots, m\}$. A subset $L$ of $[m]^n$ is a combinatorial line if there is a non-empty set $I \subseteq [n]$ and $a_i \in [m]$ for each $i \notin I$ such that 
\[L = \{(x_1, \dots, x_n) \in [m]^n : x_i = a_i  \mbox{ for all } i \notin I \mbox{ and } x_i = x_ j \mbox{ for all } i, j \in I\}.\]
The Hales--Jewett theorem is then as follows.

\begin{HJ}
For any positive integers $m$ and $r$, there exists a positive integer $n$ such that any $r$-colouring of $[m]^n$ contains a monochromatic combinatorial line.
\end{HJ}

If we define $HJ(m, r)$ to be the smallest $n$ such that the Hales--Jewett theorem holds, then the original proof results in bounds of Ackermann type for $HJ(m, r)$. In the late eighties, Shelah~\cite{S88} made a major breakthrough by finding a new way to prove the theorem which yielded primitive recursive bounds. This also gave the first primitive recursive bounds for van der Waerden's theorem. In this special case, Shelah's bound has since been drastically improved by Gowers~\cite{G01}. 

The main result of this note is a reasonable bound for the $m = 3$ case.

\begin{theorem} \label{main}
There exists a constant $c$ such that
$$HJ(3,r) \leq 2^{2^{cr}}.$$
\end{theorem}

After proving this theorem, we found that a result of this type was claimed at the end of Shelah's seminal paper. However, the brief sketch given there is at best incomplete and the bound for $HJ(3,r)$ is usually stated as being of tower type in $r$. For instance, this is the case in a paper of Graham and Solymosi~\cite{GS06} where they prove a result comparable to Theorem~\ref{main} for the coloured version of the corners theorem, that is, for finding monochromatic $(x,y)$, $(x + d, y)$, $(x, y+d)$ in any $r$-colouring of $[n]^2$. Their result is now a simple corollary of our own. We proceed straight to the details. The main idea, for those who know Shelah's proof, is to use a one-sided version of his cube lemma.


\section{The proof}

For $j = 1, \dots, t$, let $n_j = r^{6^{t-j}}$ and $s_j = n_1 + \dots + n_j$. Suppose now that $n = s_t = n_1 + \dots + n_t$ and $\chi$ is an $r$-colouring of $[3]^n$. We will show by induction, starting at $j = t$ and working downwards to $j = 0$, that there are functions $f_k : [3]^{s_j} \rightarrow \binom{0 \cup [n_k]}{2}$ for all $k > j$ such that if $w$ is a word of length $s_j$ and $f_k(w) = \{p_{k,1}, p_{k,2}\}$ with $p_{k,1} < p_{k,2}$ for each $k > j$, then the following holds:

\vspace{2mm}

For any $\ell > j$ and elements $a_{\ell + 1}, \dots, a_t$ of $[3]$, the two words $v_i = v_i(\ell; a_{\ell + 1}, \dots, a_t)$, $i = 1, 2$, have the same colour, where 
\[v_i = w \underbrace{11\dots1}_{p_{j+1, 2}}\underbrace{22\dots2}_{n_{j+1} - p_{j+1,2}} \dots \underbrace{11\dots1}_{p_{\ell, i}}\underbrace{22\dots2}_{n_{\ell} - p_{\ell,i}} \dots \underbrace{11\dots1}_{p_{t, 1}}\underbrace{a_t a_t \dots a_t}_{p_{t,2}-p_{t,1}}\underbrace{22\dots2}_{n_t - p_{t,2}}.\]
More precisely, $v_i$ is equal to $w$ for the first $s_j$ letters; for $j < k < \ell$, $v_i$ has $p_{k,2}$ ones followed by $n_k - p_{k,2}$ twos in the interval $[s_{k-1} + 1, s_k]$; in $[s_{\ell - 1} + 1, s_\ell]$, $v_i$ has $p_{\ell,i}$ ones followed by $n_\ell - p_{\ell,i}$ twos (this is the only use of the variable $i$); and, for $k > \ell$, the interval $[s_{k-1} + 1, s_k]$ consists of $p_{k, 1}$ ones, followed by $p_{k,2}-p_{k,1}$ copies of $a_k$, then by $n_k - p_{k,2}$ twos.

\vspace{2mm}

When $j = t$, there is nothing to prove. Suppose now that for any word $w'$ of length $s_{j+1}$, we have defined $f_{j+2}(w'), \dots, f_t(w')$. Let $w$ be a word of length $s_j$ and, for each $0 \leq q \leq n_{j+1}$, consider the word
\[w(q) = w\underbrace{11\dots1}_{q}\underbrace{22\dots2}_{n_{j+1} - q}\]
and write $f_k(w(q)) := \{p_{k,1}(q), p_{k, 2}(q)\}$ for all $k > j+1$. Then, for any $a_{j+2}, \dots, a_t \in [3]$, let
\[w(q; a_{j+2}, \dots, a_t) = w(q) \underbrace{11\dots1}_{p_{j+2, 1}(q)}\underbrace{a_{j+2} a_{j+2} \dots a_{j+2}}_{p_{j+2,2}(q)-p_{j+2,1}(q)}\underbrace{22\dots2}_{n_{j+2} - p_{j+2,2}(q)} \dots \underbrace{11\dots1}_{p_{t, 1}(q)}\underbrace{a_t a_t \dots a_t}_{p_{t,2}(q)-p_{t,1}(q)}\underbrace{22\dots2}_{n_t - p_{t,2}(q)}.\]
That is, $w(q; a_{j+2}, \dots, a_t)$ equals $w(q)$ for the first $s_{j+1}$ letters, then, for each $k > j+1$, the interval $[s_{k-1} + 1, s_k]$ has $p_{k,1}(q)$ ones, followed by $p_{k,2}(q) - p_{k,1}(q)$ copies of $a_k$, then by $n_k - p_{k,2}(q)$ twos.

Now, to each $w(q)$, we assign a colour $\chi_{j+1}(w(q))$, namely,
\[\prod_{k = j+2}^t f_k(w(q)) \times \prod_{a_{j+2}, \dots, a_t \in [3]} \chi(w(q; a_{j+2}, \dots, a_t)).\]
Note that the number of colours is 
\[\prod_{k=j+2}^t \binom{n_k + 1}{2} \times r^{3^{t - j - 1}} \leq (n_{j+2} \cdots n_t)^{2} r^{3^{t-j-1}}.\]
Therefore, since $n_{j+1} \geq (n_{j+2} \cdots n_t)^{2} r^{3^{t-j-1}}$, we see that there must exist two choices $p_{j+1, 1}$ and $p_{j+2, 2}$ for $q$ with $p_{j+1, 1} < p_{j+2, 2}$ such that $\chi_{j+1}(w(p_{j+1, 1})) = \chi_{j+1}(w(p_{j+1, 2}))$.

We now claim that letting $f_{j+1}(w) = \{p_{j+1, 1}, p_{j+2, 2}\}$ and, for $k > j+1$, $f_k(w) = f_k(w(p_{j+1,1})) = f_k(w(p_{j+1,2}))$ suffices. For $\ell = j + 1$, the $v_i(j+1; a_{j+2}, \dots, a_t)$ receive the same colour by the choice of $p_{j+1, 1}$ and $p_{j+1, 2}$. For $\ell > j + 1$, first note that $v_i(\ell; a_{\ell+1}, \dots, a_t)$ is the same as $v'_i(\ell; a_{\ell+1}, \dots, a_t)$, where $v'_i$ is defined relative to the word $w(p_{j+1,2})$. But, by induction, the $v'_i(\ell; a_{\ell+1}, \dots, a_t)$ receive the same colour for all choices of $\ell > j + 1$ and all $a_{\ell+1}, \dots, a_t \in [3]$. This completes our induction.

Continuing our induction all the way to $j = 0$ gives, for all $k = 1, \dots, t$, numbers $p_{k,1}, p_{k,2}$ with $0 \leq p_{k,1} < p_{k,2} \leq n_k$ such that

\vspace{2mm}

For any $0 \leq \ell \leq t$ and elements $a_{\ell + 1}, \dots, a_t$ of $[3]$, the two words $v_i = v_i(\ell; a_{\ell + 1}, \dots, a_t)$, $i = 1, 2$, have the same colour, where 
\[v_i = \underbrace{11\dots1}_{p_{1, 2}}\underbrace{22\dots2}_{n_{1} - p_{1,2}} \dots \underbrace{11\dots1}_{p_{\ell, i}}\underbrace{22\dots2}_{n_{\ell} - p_{\ell,i}} \dots \underbrace{11\dots1}_{p_{t, 1}}\underbrace{a_t a_t \dots a_t}_{p_{t,2}-p_{t,1}}\underbrace{22\dots2}_{n_t - p_{t,2}}.\]
In words, for $1 \leq k < \ell$, $v_i$ has $p_{k,2}$ ones followed by $n_k - p_{k,2}$ twos in the interval $[s_{k-1} + 1, s_k]$; in $[s_{\ell - 1} + 1, s_\ell]$, $v_i$ has $p_{\ell,i}$ ones followed by $n_\ell - p_{\ell,i}$ twos; and, for $k > \ell$, the interval $[s_{k-1} + 1, s_k]$ consists of $p_{k, 1}$ ones, followed by $p_{k,2}-p_{k,1}$ copies of $a_k$, then by $n_k - p_{k,2}$ twos.

\vspace{2mm}

To conclude the proof, take $t = r$ and consider the words $v(q) = v(0; 1, 1, \dots, 1, 3, 3, \dots, 3)$ (the $i$ is redundant when $\ell = 0$) where there are $q$ ones followed by $r - q$ threes for some $0 \leq q \leq r$. By the pigeonhole principle, there must exist $q_1$ and $q_2$ with $q_1 < q_2$ such that $\chi(v(q_1)) = \chi(v(q_2))$. But then
\begin{align*}
\chi(v(q_2)) & = \chi(v_2(q_2;3,3, \dots, 3))\\
& = \chi(v_1(q_2; 3, 3, \dots, 3))\\
& = \chi(v_2(q_2 - 1; 2, 3, 3, \dots, 3))\\
& = \dots \\
& = \chi(v_2(q_1 + 1; 2, 2, \dots, 2, 3, 3, \dots, 3)\\
& = \chi(v_1(q_1 + 1; 2, 2, \dots, 2, 3, 3, \dots, 3)
\end{align*}
where there are always $r - q_2$ threes. Here, alternate lines follow from identification between words and from an application of the conclusion above. But
\[v_1(q_1 + 1; 2, 2, \dots, 2, 3, 3, \dots, 3) = v(0; 1, 1, \dots, 1, 2, 2, \dots, 2, 3, 3, \dots, 3),\]
where there are $q_1$ ones, followed by $q_2 - q_1$ twos, then $r- q_2$ threes, so, together with $v(q_1)$ and $v(q_2)$, where the twos are replaced with threes and ones, respectively, we get the required monochromatic combinatorial line.

\vspace{1mm}

As a closing remark, we note that a similar iteration with Theorem~\ref{main} as a base shows that $HJ(4,r)$ is at most a tower of twos of height $O(r)$, improving also the bound for lines of length four.

\end{document}